\documentclass[12 pt]{article}
\usepackage{amsfonts}
\usepackage[pctex32]{graphics}
\title{Partial Covering Arrays and a Generalized Erd\H os-Ko-Rado Property}
\author{Patricia A.~Carey and Anant P.~Godbole\\
Department of Mathematics\\
East Tennessee State University}
\begin{document}
\def\qed{\vbox{\hrule\hbox{\vrule\kern3pt\vbox{\kern6pt}\kern3pt\vrule}\hrule}}
\def\ep{\varepsilon}
\def\lr{\left(}
\def\lf{\lfloor}
\def\rf{\rfloor}
\def\lc{\left\{}
\def\rc{\right\}}
\def\rr{\right)}
\def\p{\mathbb P}
\def\v{\mathbb V}
\def\a{\alpha}
\def\e{\mathbb E}
\def\l{\mathbb L}
\def\lg{{\rm lg}}
\providecommand{\floor}[1]{\left\lfloor#1\right\rfloor}
\newtheorem{thm}{Theorem}
\newtheorem{lemma}[thm]{Lemma}
\newtheorem{prop}[thm]{Proposition}
\maketitle
\begin{abstract}  The classical Erd\H os-Ko-Rado theorem states that if $k\le\floor{n/2}$ then the largest family of pairwise intersecting $k$-subsets of $[n]=\{0,1,\ldots,n\}$ is of size ${{n-1}\choose{k-1}}$.  A family of $k$ subsets satisfying this pairwise intersecting property is called an EKR family.  We generalize the EKR property and provide asymptotic lower bounds on the size of the largest family ${\cal A}$ of $k$-subsets of $[n]$  that satisfies the following property:  For each $A,B,C\in{\cal A}$, each of the four sets $A\cap B\cap C;A\cap B\cap C^C; A\cap B^C\cap C; A^C\cap B\cap C$ are non-empty.    This generalized EKR (GEKR) property is motivated, generalizations are suggested, and a comparison is made with fixed weight $3$-covering arrays.  Our techniques are probabilistic, and reminiscent of those used in \cite{gss} and in the work of Roux, as cited in \cite{sloane}.  
\end{abstract}
\section{Introduction}  The classical Erd\H os-Ko-Rado (EKR) theorem (\cite{vanlint}) states that if $k\le\floor{n/2}$ then the largest family of pairwise intersecting $k$-subsets of $[n]=\{0,1,\ldots,n\}$ is of size ${{n-1}\choose{k-1}}$, and is given, moreover, by the class of subsets of $[n]$ containing a fixed element $a$.  A family of pairwise intersecting $k$ subsets is called an EKR family.  We generalize the EKR property and provide asymptotic lower bounds on the size of the largest family ${\cal A}$ of $k$-subsets of $[n]$  that satisfies the following property:  For each $A,B,C\in{\cal A}$, each of the four sets $A\cap B\cap C;A\cap B\cap C^C; A\cap B^C\cap C; A^C\cap B\cap C$ are non-empty, where $A^C$ denotes the complement of the set $A$.  

Now why might such a property be of interest?  Here is motivation for our choice of this somewhat unusual generalized EKR (GEKR) property, together with possible extensions.  The EKR theorem can be thought of in several ways, but one is the following:  If each person speaks a different set of $k$ languages out of a total of $n$, what is the largest number of people that can have two-way conversations with each other?  Or, in a ``small world network" context, if each person knows a unique set of $k$ others, what is the largest number of people possible so that any two have a mutual acquaintance?  Now imagine that we desire a situation in which any two of any three people can have a conversation with no fear of the third eavesdropping, and yet all three are able to communicate if necessary.  Or one in which any two of any three people have a common friend who is a stranger to the third, while there also exists a person who is a mutual friend of all three.  We seek therefore, to construct a family of ``partially 3-independent sets" (see, e.g., \cite{anderson} for the definition of $k$-independence.) Our definition can be extended to one in which, given any three people, $1\le\beta\le 8$ specified regions in the associated Venn diagram are non-empty.   This situation can clearly be generalized to more than three sets, but we choose not to do so. Another possibility might be (using the language analogy) the following:  Language ability is at three levels: ignorant (0); novice (1), and expert (2).  Two persons who are at the novice and expert level at a particular language can converse at the novice level, and thus keep keep their conversation  secret from the third, assumed to be ignorant.   Likewise, two persons at the level `2' can speak rapidly and at a high level, thus keeping their conversation undecipherable to a third who is at the `0' or `1' level.  Three persons can converse if there exists a language at which each is at level 1 or above.  

Key to our development is a basic comparison between the GEKR property and fixed weight binary 3-covering arrays (\cite{colbourn} and \cite{sloane} are comprehensive surveys, and both provide an exhaustive list of references on covering arrays).  Recall that a $t$-covering array ${\cal C}(n,q, t,\lambda)$ is defined as an $m\times n$ array satisfying the property that for any choice of $t$ rows, each of the $q^t$ $q$-ary $t$-tuples appear at least $\lambda$ times among the columns of the selected rows.  We define the partial $t$-covering array ${\cal C}(n,q, t, A\subseteq\{0,1,\ldots,q-1\}^t,\lambda)$ to be an $m\times n$ array satisfying the property that for any choice of $t$ rows, each of the $q$-ary $t$-tuples in the ensemble $A$ appear at least $\lambda$ times among the columns of the selected rows.  We thus see that the GEKR property holds for $m$ ``people" if and only if the $m\times n$ ``person-language" incidence matrix is a partial ${\cal C}(n,2, 3, A, 1)$ covering array with $A=\{(0,1,1)\cup(1,0,1)\cup (1,1,0)\cup (1,1,1))\}$, and with each row having weight $k$.  We seek to find lower bounds on the maximum number of rows $m$ so that such an array exists, and our results are typically of the form $m\ge C^n$ for some constant $C$.
Our techniques are probabilistic, and reminiscent of those used in \cite{gss} and by Roux, the latter as cited in \cite{sloane}.  Connections between the EKR property and $t$-covering arrays have been exploited in a different context by researchers such as Karen Meagher, Lucia Moura, and Brett Stevens in Ottawa.  For example, in \cite{meagher2}, EKR theorems for are proved for set partitions. These are connected to strength-2 covering arrays with any alphabet rather than strength-3 binary covering arrays.  See also \cite{meagher}.

Our results are presented in the next section.  Several other generalizations of the EKR theorem are studied in \cite{anderson}, \cite{erdos}, \cite{meagher2}, and \cite{uiuc}; the references contained in these sources span, between them, almost 45 years of developments. 
\section{Results}  We start by considering a result in which each row of the $m\times n$ array is allowed to have {\it expected} weight $k$.  Specifically, the strategy is to independently place a one in any of the $m\cdot n$ places with probability $\alpha=k/n$ and a zero with probability $1-\alpha$.  A set of three rows is said to be {\it deficient} if the rows do not contain one of the four vectors $(0,1,1),(1,0,1),(1,1,0), (1,1,1)$ among their columns. Set $X=\sum_{j=1}^{m\choose 3}I_j$, where $I_j$ equals one or zero according as the $j$th set of three rows is deficient (or not).  It is evident that $\{X=0\}$ if and only if the array satisfies the GEKR property, so that $\p(X=0)>0\Leftrightarrow$ the appropriate partial 3-covering array can be constructed. The Lov\'asz local lemma is used extensively throughout the paper -- it provides a simple
 condition which guarantees a positive probability for the event that no set of three rows is deficient; see \cite{alon} for this result, stated below for convenience:
\begin{lemma}
Let $A_1$, $A_2$,\ldots, $A_N$ be events in an arbitrary probability
space.  Suppose that each event $A_i$ is mutually 
independent of a set of all the other events $A_j$ but at most $d$, and
that $\p(A_i) \leq p$ for all $1 \leq i \leq N$.
 If $ep(d+1) \leq 1$ then $\p(\bigcap_{i=1}^N A_i^C) > 0$. 
\end{lemma}
\begin{thm}  Consider an $m\times n$ array of zeros and ones with each entry being independently chosen to be a one with probability $\alpha$ and a zero with probability $1-\a$.  Let $X$ be the number of sets of three deficient rows.  Let $A=\{(0,1,1)\cup(1,0,1)\cup (1,1,0)\cup (1,1,1)\}$.  Then,  $$m\le{\sqrt{\frac{2}{3e}}}\lr\frac{1}{p_n(\a)}\rr^{\frac{1}{2}}=\zeta(n)\Rightarrow \p(X=0)>0,$$ where
$p_n(\a)=(1-\a^3)^n+3(1-\a^2(1-\a))^n$, so that the largest number of rows in a  ${\cal C}(n,2, 3, A, 1)$ covering array, with each row having expected weight $k=\a\cdot n$, is at least $\zeta(n)$.
\end{thm}

\medskip

\noindent{\bf Proof}  We use the Lov\'asz local lemma.  Let $A_i$ be the event that the $i$th set of 3 rows is deficient, so that we have, for each $i$, 
\begin{eqnarray*}
\p(A_i)&=&\p\lr\bigcup_{j\in A}\{j\ {\rm is\ missing}\}\rr\\
&\le&\p((1,1,1)\ {\rm is\ missing})+3\p((1,1,0)\ {\rm is\ missing})\\
&=&(1-\a^3)^n+3(1-\a^2(1-\a))^n.\\
\end{eqnarray*}
Now, it is evident that the dependence number $d$ satisfies $d+1\le3{{m-1}\choose{2}}\le 3m^2/2$, so that $\p(X=0)>0$ provided that
$${{3e}\over{2}}m^2p_n(\a)\le1,$$ i.e., if 
$$m\le {\sqrt{\frac{2}{3e}}}\lr\frac{1}{p_n(\a)}\rr^{\frac{1}{2}},$$
as asserted.\hfill\qed  

\medskip

\noindent{\bf Remarks.} Notice that
$$p_n(\a)=\cases{4\lr\frac{7}{8}\rr^n&if\ $\a=\frac{1}{2}$\cr
3(1-\a^2(1-\a))^n(1+o(1))&if\ $\a>\frac{1}{2}$\cr
(1-\a^3)^n(1+o(1))&if\ $\a<\frac{1}{2}$.\cr}$$  It follows, due to the fact that $(1-\a^3)^n$ is monotone decreasing and $3(1-\a^2(1-\a))^n$ is decreasing on the interval $[1/2,2/3]$, that the function $p_n(\a)$ attains its minimum at $\a=2/3$ -- yielding the conclusion that the maximum size of a ${\cal C}(n,2, 3, A, 1)$ covering array appears to be when the expected weight of each row is $\frac{2}{3}\cdot n$.  Also, a lower bound on this size is 
$${\sqrt{\frac{2}{9e}}}\cdot(1.0834\ldots)^n.$$  Furthermore, the rate of growth of the lower bound is exponential no matter what $\a>0$ is.  It is also worth mentioning that nowhere in this paper can one obtain  improved {\it asymptotic} results by using a more careful estimation of the probability of sets of rows being deficient, using, for example, the inclusion-exclusion principle.

Might we have been able to incorporate into our analysis the rather realistic case where $\a=\a_n\to\infty; \a_n/n\to0$?  Since in this case $p_n^{1/2}(\a)\sim(1-\a_n^3)^{n/2}\sim\exp\{-n\a_n^3/2\},$ we need, at the very least, to have $n\a_n^3\to\infty$ for the size of the array to grow to infinity with $n$.  But the results are entirely satisfactory in this case.  For example, with $n=10^9$ and $\a_n=\log^{2/3}n/n^{1/3}=.0075$, we get a lower bound of $\sim2\cdot10^{91}$ on the size of the array.

Next we turn to our main result, which, while being similar in spirit to Theorem 2, is based on the GEKR property holding for {\it fixed weight} rows.  The analysis gets more complicated due to the fact that entries within each row are no longer independent, but the dependence structure between rows stays the same.  We shall see, moreover, that a better bound is obtained in the fixed weight case unless both $\a$ and $n$ are small.
\begin{thm}
Consider an $m\times n$ array of zeros and ones with each row containing $k=\a n$ ones and $n-k=(1-\a)n$ zeros $(\a\in{\mathbb Q}\cap(0,1); \a n\in{\mathbb Z}^+)$, so that each of the ${n\choose k}$ configurations of ones and zeros are equally likely. Let $X$ be the number of sets of three deficient rows.  Let $A=\{(0,1,1)\cup(1,0,1)\cup (1,1,0)\cup (1,1,1)\}$.  Then there exists $K>0$ such that  $$m\le\lr\frac{1}{\sqrt{K}\mu(\a)^{n/2}n^{1/4}}\rr(1+o(1)):=\nu_n(\a)\Rightarrow \p(X=0)>0,$$ where
$\mu(\a)$ is given on the last line of the proof that follows.  In other words, the largest number of rows in a  ${\cal C}(n,2, 3, A, 1)$ covering array is at least $\nu_n$.
\end{thm}

\medskip

\noindent{\bf Proof}  We set $X=\sum_{j=1}^{m\choose 3}I_j$; as before, $X$ is the number of deficient sets of three rows.  As before, the $d$ in the Lov\'asz local lemma can be bounded by $3m^2/2$, and the main problem is to estimate the probability that a given set of rows is deficient.  We have
\begin{eqnarray*}p&=&\p(I_j=1)\le\p((1,1,1) {\rm is\ missing}) + 3\p((1,1,0) {\rm is\ missing})\nonumber\\
&=&\sum_u\frac{{n\choose{\a n}}{{\a n}\choose{u}}{{n-\a n}\choose{\a n-u}}{{n-u}\choose{\a n}}}{{n\choose{\a n}}^3}+3\sum_u\frac{{n\choose{\a n}}{{\a n}\choose{u}}{{n-\a n}\choose{\a n-u}}{{n-u}\choose{n-\a n}}}{{n\choose{\a n}}^3}\nonumber\\
&=&\sum_u\frac{{{r}\choose{u}}{{n-r}\choose{r-u}}{{n-u}\choose{r}}}{{n\choose{r}}^2}+3\sum_u\frac{{{r}\choose{u}}{{n-r}\choose{r-u}}{{n-u}\choose{n-r}}}{{n\choose{r}}^2}\qquad(\a n=r)\nonumber\\
&=&\sum_1\phi(u)+\sum_2\psi(u)\quad{\rm(say).}
\end{eqnarray*}
Consider $\Sigma_1$ first.  We must have, in this first sum,
\[
u\le\min\{\a n, (1-\a)n\}=\min\{r,n-r\};
\]
\[u\ge\max\{0,(2\a-1)n\}=\max\{0,2r-n\}.
\]
With this in mind, consider the ratio
\[\frac{\phi(u+1)}{\phi(u)}=\frac{{r \choose u+1}{n-r\choose r-u-1}{n-u-1\choose r}}{{r\choose u}{n-r\choose r-u}{n-u\choose r}}\]
of successive terms, which may be checked to exceed one if and only if
\[(n-r+2)u^2+(r^2-2r-n^2-n+1)u+(nr^2-r^3-n^2+2nr-n)\ge0.\]
The idea is to find which value of $u$ maximizes $\phi(u)$.  Towards this end, we see that $\phi$ is increasing if and only if
$u\le u_1; u\ge u_2,$
where the smaller and larger roots $u_1,u_2$ of the above quadratic are respectively given by
\[u_1=\frac{1}{2(n-r+2)}\cdot\bigg(-r^2+2r+n^2+n-1
-{\sqrt{\gamma(n,r)}}\bigg),\]
and 
\[u_2=\frac{1}{2(n-r+2)}\cdot\bigg(-r^2+2r+n^2+n-1
+{\sqrt{\gamma(n,r)}}\bigg),\]
with 
\begin{eqnarray*}\gamma(n,r)&=&1-2nr^2-16nr-4r+6n+6r^2+11n^2\\
{}&&+4r^3-6r^2n^2-3r^4+n^4+6n^3-8rn^2+8r^3n.\end{eqnarray*}
We simplify by reintroducing the parameter $\a$; $\a n=r$, and setting
\[e=1-3\a^4-6\a^2+8\a^3;\]
\[f=-2\a^2+4\a^3+6-8\a;\]
$e$ and $f$ are, respectively, the coefficients of the fourth and third order terms in the radical above.  We thus get, noting that $e>0$ for $\a<1$,
\[u_1=\frac{1}{(2n-2\a n+4)}\cdot\bigg(n^2(1-\a^2)+n(1+2\a)-1
-{\sqrt{e}}n^2-\frac{f}{2\sqrt{e}}n+O(1)\bigg)\]
and
\[u_2=\frac{1}{(2n-2\a n+4)}\cdot\bigg(n^2(1-\a^2)+n(1+2\a)-1
+{\sqrt{e}}n^2+\frac{f}{2\sqrt{e}}n+O(1)\bigg).\]
The critical points we need to investigate, together with the endpoints of the summation, are thus (within $O(1)$ of)
\begin{eqnarray*}u_{1,{\rm approx}}&=&\frac{(1-\a^2-\sqrt{e})}{(2-2\a )}\cdot n\\&=&\frac{(1-\a^2-\sqrt{1-3\a^4-6\a^2+8\a^3})}{(2-2\a )}\cdot n\\
&:=&\beta\cdot n=\beta(\a)n\end{eqnarray*}
and
\[u_{2,{\rm approx}}=\frac{(1-\a^2+\sqrt{e})}{(2-2\a )}\cdot n=\frac{(1-\a^2+\sqrt{1-3\a^4-6\a^2+8\a^3})}{(2-2\a )}\cdot n.\]
Consider Figure 1 below which graphs respectively $1/n$ times (i) the lower limit of summation of $\Sigma_1$, given by the straight line; (ii) the upper limit of summation of $\Sigma_1$, given by the triangular plot; and (iii) $u_{1,{\rm approx}}$, given by the increasing curve; and (iv) $u_{2,{\rm approx}}$, the unimodal curve. Note the interesting (but obvious in hindsight) fact that the lower limit of summation exceeds the upper limit if $\a>2/3$ and thus $\p((1,1,1)\ {\rm is\ missing}) $ equals zero in this case.
Thus the maximum of $\phi(u)$ is attained at $u_{1,{\rm approx}}+O(1)$ for each $\a\le2/3$.

\bigskip
\hfil\scalebox{1}{\includegraphics{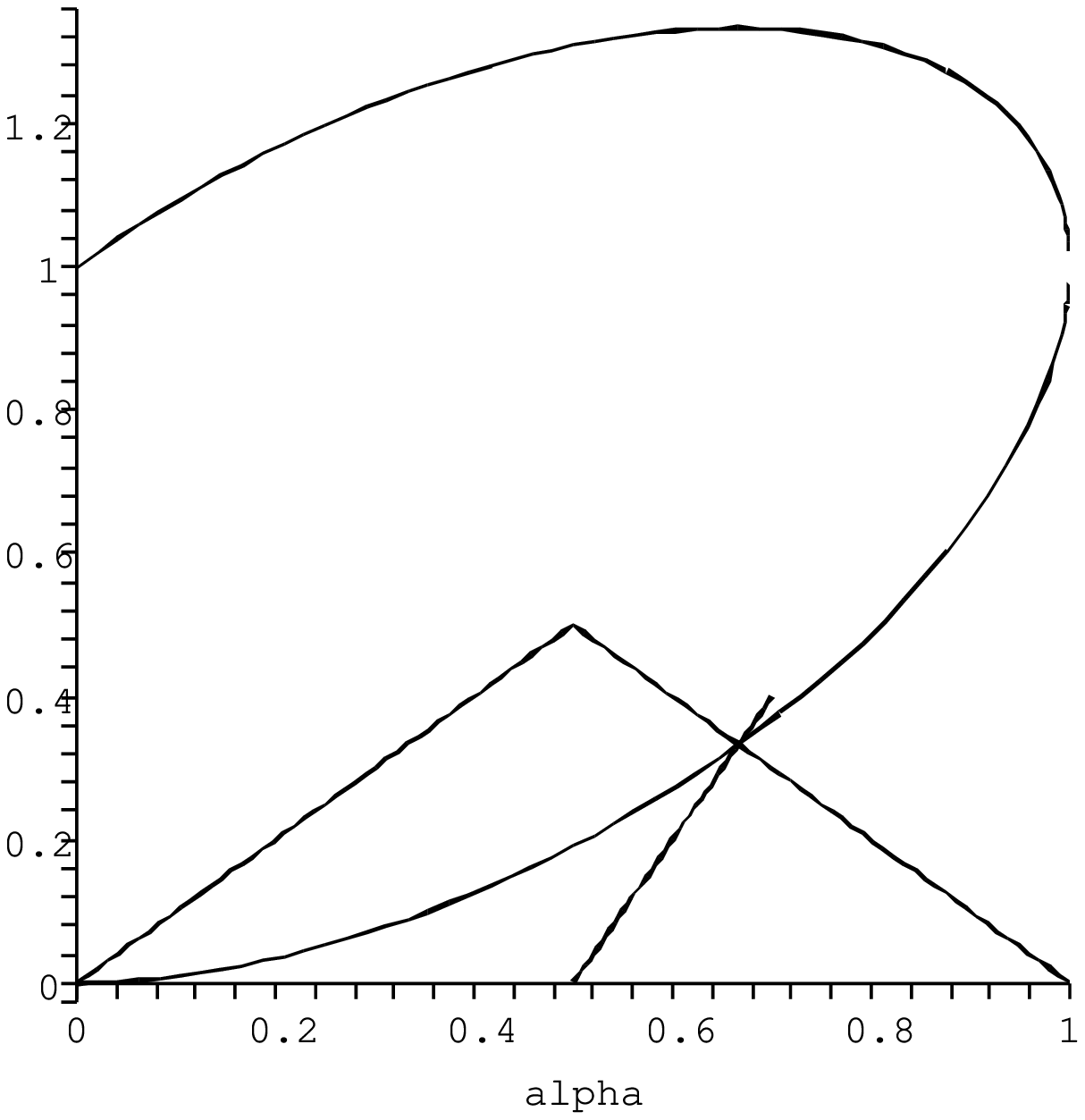}}\hfil
\centerline{Figure 1}

\bigskip

Consider $\Sigma_2$, in which the limits of summation are
\[
0\vee(2\a-1)n\le u\le\a n.
\]
Now we find that 
$$\frac{\psi(u+1)}{\psi(u)}=\frac{(\a n-u)^3}{(u+1)(n-2\a n+u+1)(n-u)}\ge1$$
if and only if $au^2+bu+c\ge0$, where
\[a=\a n+2;\]
\[b=-3\a^2n^2-n^2+2\a n^2-n-2\a n+1;\]
and
\[c=\a^3n^3-n^2+2\a n^2-n.\]
The roots $v_1$ and $v_2$ of the quadratic $au^2+bu+c=0$ are
\[\frac{1}{2(\a n+2)}\lr3\a^2n^2+n^2-2\a n^2+n+2\a n-1\pm\sqrt{\delta(n,\a)}\rr\]
where
\begin{eqnarray*}\delta(n,\a)&=&1-4\a n-2\a^2n^2-4\a n^2+4\a^3n^3+7n^2+6n+5\a^4n^4\\
&&{}+10\a^2n^4-12\a^3n^4-10\a^2n^3-4\a n^4+2n^3+n^4+4\a n^3.\end{eqnarray*}
We thus see as before that $v_1$ and $v_2$ can be approximated by 
\begin{eqnarray*}v_{1,{\rm approx}}&=&\frac{3\a^2n+n-2\a n-\sqrt{1+10\a^2-12\a^3-4\a+5\a^4}n}{2\a}\\&:=&\kappa n=\kappa(\a)n\end{eqnarray*}
and
\[v_{2,{\rm approx}}=\frac{3\a^2n+n-2\a n+\sqrt{1+10\a^2-12\a^3-4\a+5\a^4}n}{2\a}.\]

\bigskip
\hfil\scalebox{1}{\includegraphics{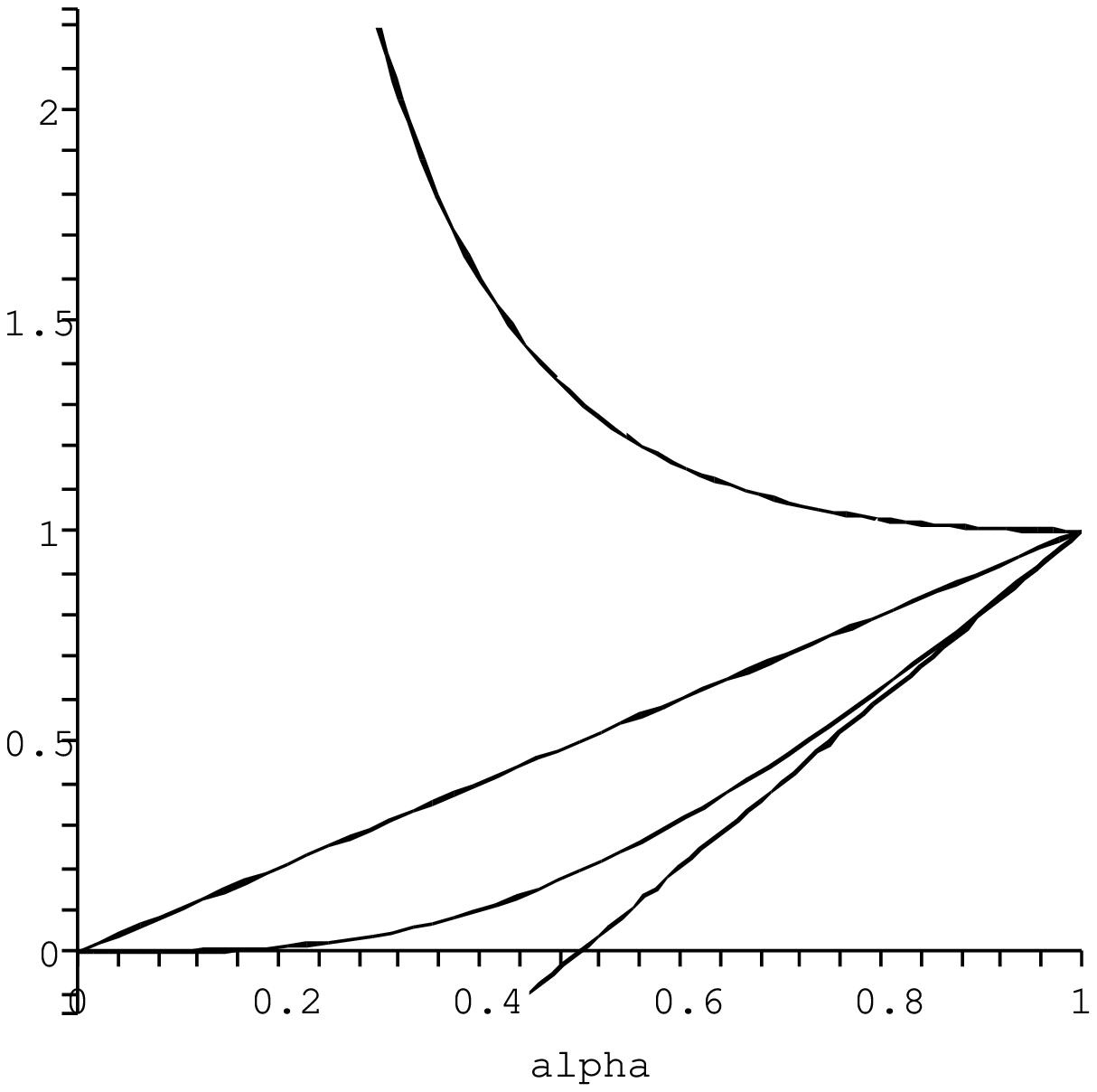}}\hfil
\centerline{Figure 2}

\bigskip
\noindent The functions $v_{1,{\rm approx}}/n; v_{2,{\rm approx}}/n; \a$ and $(2\a-1)$ are plotted in Figure 2.  These graphs reveal that the maximum of $\psi(u)$ occurs at $v_{1,{\rm approx}}+O(1)$.

Returning to $\Sigma_1$, we estimate as follows; throughout this paper $K$ will denote a generic constant whose value might change from line to line.
\begin{eqnarray*}
\sum_1\phi(u)&\le&Kn\cdot\max\phi(u)\\
&\le&Kn\cdot\phi(u_{1,{\rm approx}}+O(1))\\
&\le&Kn{{{{\a n}\choose{\beta n}}{{n-\a n}\choose{\a n-\beta n}}{{n-\beta n}\choose{\a n}}}\over{{n}\choose {\a n}}^2}.
\end{eqnarray*}
Stirling's approximation $N!\approx\sqrt{2\pi N}(N/e)^N$, applied to each of the binomial coefficients above, yields after some simplification, 
\begin{eqnarray*}
\sum_1\phi(u)&\le&Kn\phi(\beta n)\\
&=&K\sqrt{n}\lr\frac{(1-\a)^{3-3\a}(1-\beta)^{1-\beta}\a^{2\a}}{\beta^{\beta}(\a-\beta)^{2\a-2\beta}(1-2\a+\beta)^{1-2\a+\beta}(1-\a-\beta)^{1-\a-\beta}}\rr^n\\
&=&K\sqrt {n}\xi^n\quad{\rm say}.
\end{eqnarray*}

We treat $\Sigma_2$ in a similar fashion: For each $\a\in[0,1]$, 
\begin{eqnarray*}
\sum_2\psi(u)&\le&Kn\cdot\max\psi(u)\\
&\le&Kn\cdot\psi(v_{1,{\rm approx}})\\
&=&{{{{\a n}\choose{\kappa n}}{{n-\a n}\choose{\a n-\kappa n}}{{n-\kappa n}\choose{n-\a n}}}\over{{n}\choose {\a n}}^2}.
\end{eqnarray*}
Stirling's approximation again yields after some simplification, 
\begin{eqnarray*}
\sum_2\psi(u)&\le&Kn\phi(\kappa n)\\
&=&K\sqrt{n}\lr\frac{\a^{3\a}(1-\a)^{2-2\a}(1-\kappa)^{1-\kappa}}{\kappa^\kappa(\a-\kappa)^{3\a-3\kappa}(1-2\a+\kappa)^{1-2\a+\kappa}}\rr^n\\
&=&K\sqrt{n}\theta^n\quad{\rm say}.
\end{eqnarray*}
In Figure 3 we have plotted $\xi(\a)$ (the concave down curve) and $\theta(\a)$; note that the domain of $\xi$ is $[0,2/3]$, while $\theta$ is defined on $[0,1]$.  It follows that the $\Sigma_1$ sum dominates the probability of three specific rows being deficient if $\a\le1/2$, with $\Sigma_2$ taking over for $\a>1/2$ {\it just as in the case of rows with expected weight $k$}. There is a significant difference in the two results, however.  Figure 4, which plots $\theta$ in the vicinity of $\a=0.74$, reveals that the minimum of $\theta$ is attained around $\a=0.7395$ in contrast to the fact that the minimum was attained at $\a=2/3$ in the ``independent" case.  Thus our proof suggests that the maximal size of a generalized EKR array occurs when $\a\approx0.7395$.

\bigskip
\hfil\scalebox{1}{\includegraphics{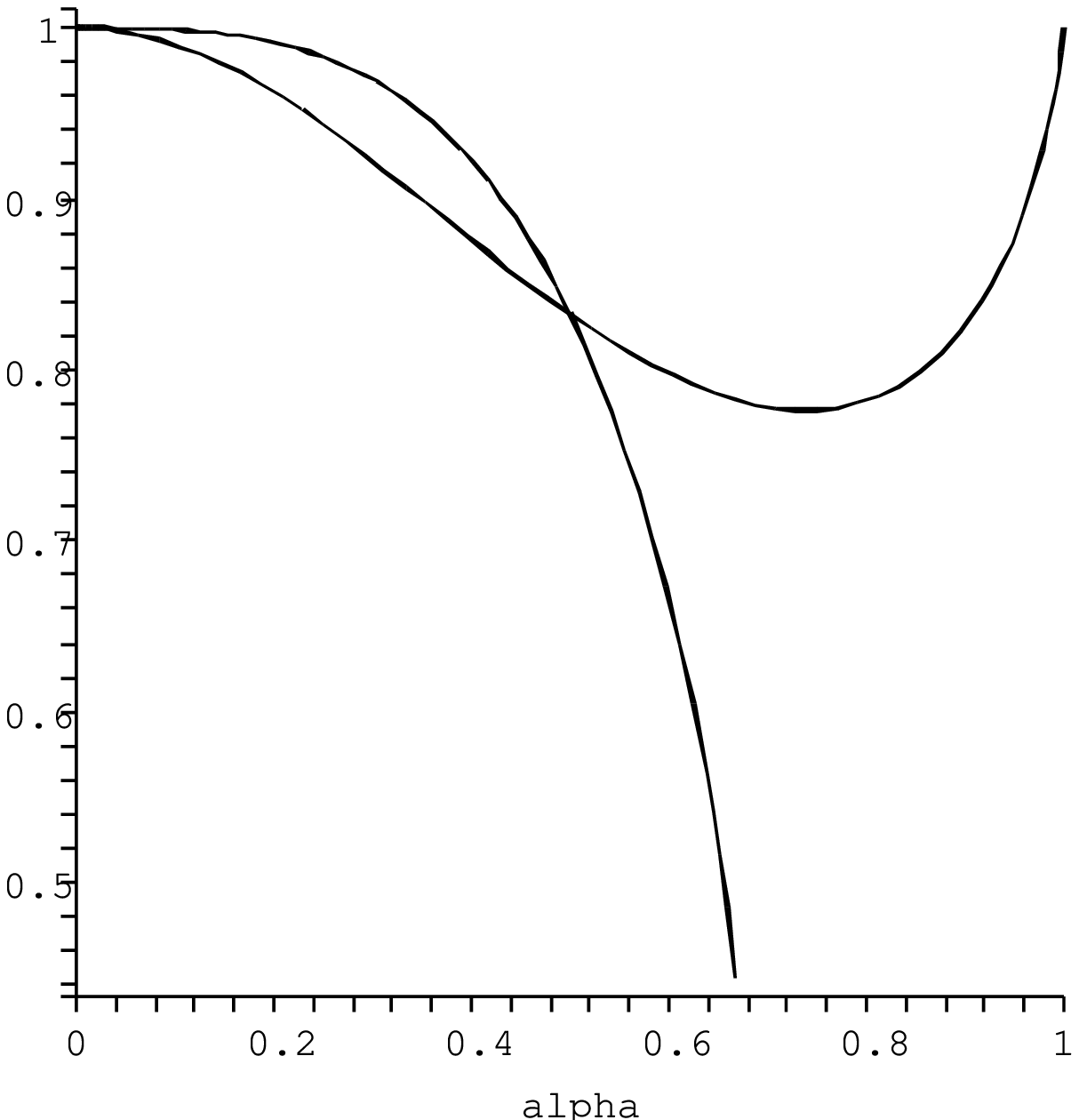}}\hfil
\centerline{Figure 3}

\bigskip
\bigskip

\hfil\scalebox{1}{\includegraphics{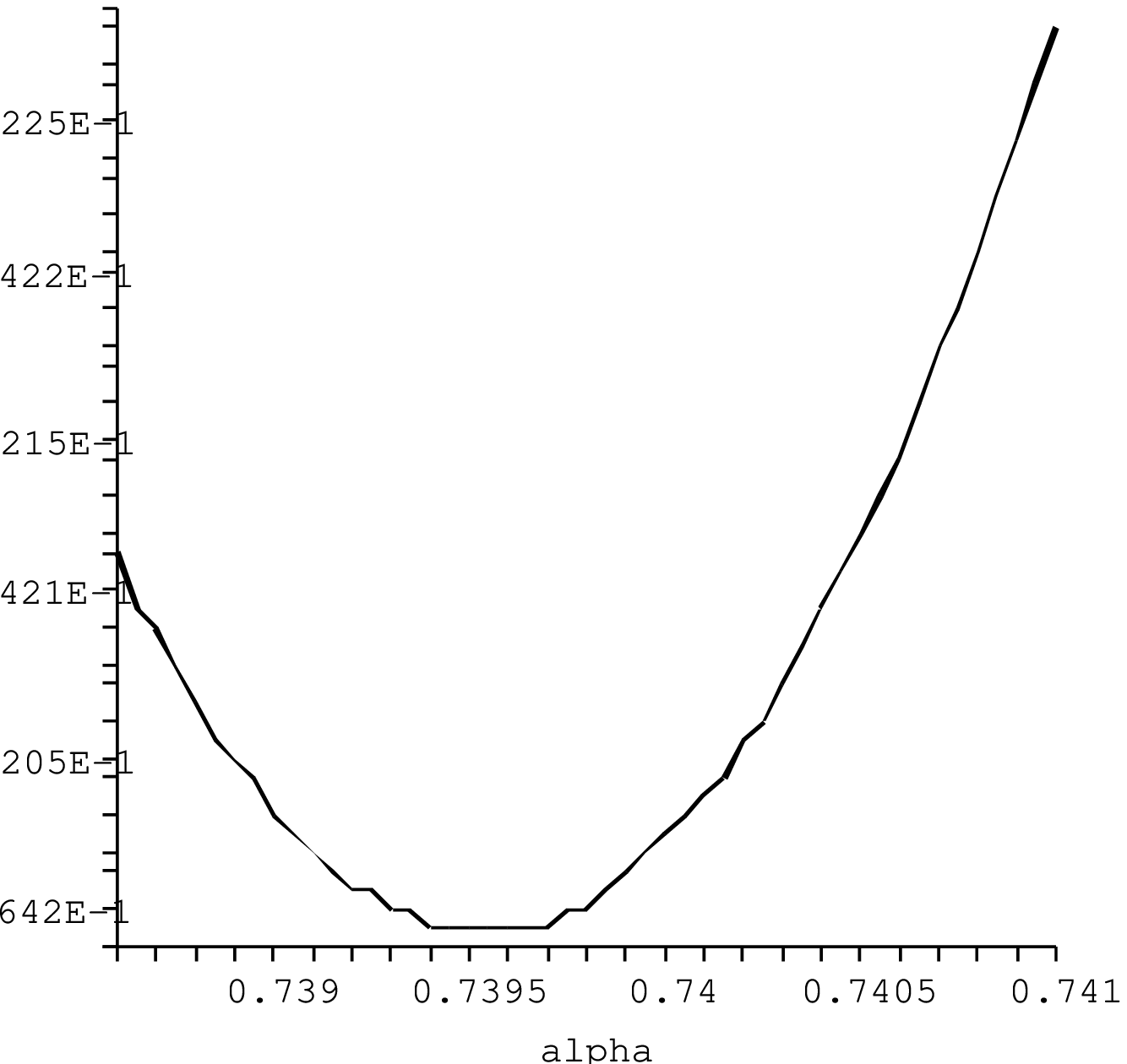}}\hfil
\centerline{Figure 4}

\bigskip

We now finish off the proof of Theorem 3.  By the Lov\'asz local lemma, $\p(X=0)>0$ whenever $3epm^2/2\le1$, or if 
\[(1+o(1))K\sqrt{n}\xi^nm^2\le1, \quad(\a\le1/2);\]
\[(1+o(1))K\sqrt{n}\theta^nm^2\le1, \quad(\a>1/2).\]
We now simply set $\mu=\xi$ if $\a\le1/2$ and $\mu=\theta$ if $\a>1/2$.  Theorem 3 follows.\hfill\qed

\medskip

\noindent{\bf A Remark and Numerical Values.}  First note that we recover the result of Roux, as discussed in \cite{sloane}, on setting $\a=1/2$ (this is the {\it highly} symmetrical case studied earlier in the literature.)  In this case we get
\[u_{1,{\rm approx}}=v_{1,{\rm approx}}=\frac{3-\sqrt {5}}{4}\cdot n,\quad{\rm i.e.}\ \beta=\kappa=\frac{3-\sqrt {5}}{4}.\]
Second, we provide below in Tables 1 and 2 some numerical values for various values of $\a$,  for the independent and fixed weight models studied in Theorems 2 and 3 respectively. While computing the bounds given by Theorem 3, we have ignored the effect of the $(1+o(1))$ term, and pretended that $K=1$.  The ``interesting" {\it first} choices for $\a$ are each roughly 1/6, correspond to $n=10,000$, and yield $m\approx 6.5\cdot10^9$ -- the world's current population.
 
\bigskip
\bigskip

 \centerline{Table 1}
 \centerline{\it The Independent Model}
 \medskip
$$\vbox{\halign{
\hfil#\hfil&\qquad
\hfil#\hfil&\qquad
\hfil#\hfil&\qquad
\hfil#\hfil&\qquad
\hfil#\hfil\cr
$\a/n$& 10,000& 100,000& 300,000& 1,000,000\cr
0.1669& $6.51\cdot10^9$& $7.66\cdot10^{100}$& $1.83\cdot10^{303}$& $3.88\cdot10^{1011}$\cr
0.2& $1.37\cdot10^{17}$& $1.29\cdot10^{174}$& $8.79\cdot10^{522}$& $7.22\cdot10^{1743}$\cr
1/3& $4.34\cdot10^{81}$& $1.64\cdot10^{819}$& $1.81\cdot10^{2458}$& $8.00\cdot10^{8194}$\cr
0.5& $2.26\cdot10^{289}$& $9.80\cdot10^{2898}$& $1.53\cdot10^{8698}$& $2.33\cdot10^{28995}$\cr
2/3& $4.32\cdot10^{347}$& $1.79\cdot10^{3481}$& $7.00\cdot10^{10444}$& $2.63\cdot10^{34817}$\cr
0.7395& $1.50\cdot10^{333}$& $4.61\cdot10^{3336}$& $1.20\cdot10^{10011}$& $3.37\cdot10^{33371}$\cr
0.8& $7.47\cdot10^{296}$& $4.28\cdot10^{2973}$& $9.63\cdot10^{8921}$& $1.64\cdot10^{29741}$\cr
}}$$

\bigskip
\bigskip

 \centerline{Table 2}
 \centerline{\it The Fixed Weight Model}
 \medskip
$$\vbox{\halign{
\hfil#\hfil&\qquad
\hfil#\hfil&\qquad
\hfil#\hfil&\qquad
\hfil#\hfil&\qquad
\hfil#\hfil\cr
$\a/n$& 10,000& 100,000& 300,000& 1,000,000\cr
0.1685& $6.50\cdot10^9$& $7.61\cdot10^{106}$& $1.06\cdot10^{323}$& $6.52\cdot10^{1079}$\cr
0.2& $2.32\cdot10^{17}$& $2.57\cdot10^{182}$& $4.08\cdot10^{549}$& $1.26\cdot10^{1835}$\cr
1/3& $9.50\cdot10^{92}$& $3.39\cdot10^{938}$& $9.36\cdot10^{2817}$& $1.99\cdot10^{9396}$\cr
0.5& $9.00\cdot10^{396}$& $1.97\cdot10^{3978}$& $1.84\cdot10^{11937}$& $8.82\cdot10^{39793}$\cr
2/3& $4.50\cdot10^{530}$& $1.93\cdot10^{5315}$& $1.73\cdot10^{15948}$& $7.27\cdot10^{53163}$\cr
0.7395& $3.28\cdot10^{548}$& $8.27\cdot10^{5493}$& $1.35\cdot10^{16484}$& $1.46\cdot10^{54950}$\cr
0.8& $1.74\cdot10^{533}$& $1.48\cdot10^{5341}$& $7.86\cdot10^{16025}$& $5.19\cdot10^{53422}$\cr
}}$$

\section{Open Questions} Of possible interest might be lower bounds on the size of partial $t$-covering arrays ${\cal C}(n,q, t, A,\lambda)$ for general values of the parameters.  More crucial, however, would be construction, algorithms, and improvements -- especially for small values of $n$ -- for the baseline case studied in this paper:  $t=3, q=2;\lambda=1$, and $A$ as specified by the GEKR property.  Also, we feel that the situation where $\a=\a_n$ varies as a function of $n$ so that $\a_n/n\to0$ needs further investigation. Last but not least, Karen Meagher has asked a deep question:  Can the methods in this paper be used to try to get a generalized {\it Sperner} property? Perhaps it might be possible to use a generalized Sperner type result to say that the largest system (not necessarily of $k$-sets) that satisfies the generalized EKR property is a system of $k$-sets.
\section{Acknowledgements} This research leading to this paper was conducted during the 2004--05 academic year, and was part of Carey's mandatory undergraduate research requirement at ETSU.  The research of her advisor, the second-named author, was supported by NSF Grant DMS-0139286.

\end{document}